\documentstyle[doublespace,11pt]{article}

    \csname @addtoreset\endcsname{figure}{section} 
    \csname @addtoreset\endcsname{table}{section} 
 
\newcounter{thanksnum} 
\def\thanksnumber#1 
{\setcounter{thanksnum}{\value{footnote}}\setcounter{footnote}{#1}%
                     \addtocounter{footnote}{-1}\footnotemark 
                     \setcounter{footnote}{\value{thanksnum}}} 
\def\newtheoremz#1{\@ifnextchar[{\@othmz{#1}}{\@nthmz{#1}}} 
 
\def\@nthmz#1#2{%
\@ifnextchar[{\@xnthmz{#1}{#2}}{\@ynthmz{#1}{#2}}} 
 
\def\@xnthmz#1#2[#3]{\expandafter\@ifdefinable\csname #1\endcsname 
{\@definecounter{#1}\@addtoreset{#1}{#3}%
\expandafter\xdef\csname the#1\endcsname{\expandafter\noexpand 
  \csname the#3\endcsname \@thmcountersepz \@thmcounterz{#1}}%
\global\@namedef{#1}{\@thmz{#1}{#2}}\global\@namedef{end#1}{\@endtheoremz}}} 
 
\def\@ynthmz#1#2{\expandafter\@ifdefinable\csname #1\endcsname 
{\@definecounter{#1}%
\expandafter\xdef\csname the#1\endcsname{\@thmcounterz{#1}}%
\global\@namedef{#1}{\@thm{#1}{#2}}\global\@namedef{end#1}{\@endtheoremz}}} 
 
\def\@othmz#1[#2]#3{\expandafter\@ifdefinable\csname #1\endcsname 
  {\global\@namedef{the#1}{\@nameuse{the#2}}%
\global\@namedef{#1}{\@thmz{#2}{#3}}%
\global\@namedef{end#1}{\@endtheoremz}}} 
 
\def\@thmz#1#2{\refstepcounter 
    {#1}\@ifnextchar[{\@ythmz{#1}{#2}}{\@xthmz{#1}{#2}}} 
 
\def\@xthmz#1#2{\@begintheoremz{#2}{\csname the#1\endcsname}\ignorespaces} 
\def\@ythmz#1#2[#3]{\@opargbegintheoremz{#2}{\csname 
       the#1\endcsname}{#3}\ignorespaces} 
 
\def\@thmcounterz#1{\noexpand\arabic{#1}} 
\def\@thmcountersepz{.} 
\def\@begintheoremz#1#2{ \trivlist \item[\hskip \labelsep{\bf #1\ #2}]} 
\def\@opargbegintheoremz#1#2#3{ \trivlist 
      \item[\hskip \labelsep{\bf #1\ #2\ (#3)}]} 
\def\@endtheoremz{\endtrivlist}

\newtheorem{theorem}{Theorem}
\newtheorem{lemma}{Lemma}
\newtheorem{definition}{Definition}

\def\defi{\stackrel{{\scriptscriptstyle \Delta}}{=}}

\def\a{\alpha}  
\def\d{\delta}

\def\w{\widehat}

\def\R{{\bf R}} 
 
\def\P{{\bf P}}

\def\L{L} 
 
\def\b{\beta} 
 
\def\g{\gamma} 
 
\def\W{{\cal W}^*}

\def\oo{\overline} 
 
\def\p{\partial} 
\def\G{\Gamma}

\def\V{{\cal V}}

\def\L{{\cal L}}

\newcommand{\be}{\begin{equation}} 
\newcommand{\ee}{\end{equation}} 
\newcommand{\bd}{\begin{displaymath}} 
\newcommand{\ed}{\end{displaymath}} 
\newcommand{\ba}{\bd\begin{array}{rl}} 
\newcommand{\ea}{\end{array}\ed} 
\font\sm=cmr10 
\date{\  }
\def\L{{\cal L}}

\def\W{{\cal W}}

\def\Q{{\cal Q}}

\def\R{{\bf R}}

\title
{
A new property of absorbed diffusions\thanks{
Supported by a University research grant,
The Hong Kong Polytechnic
University, and  Russian Foundation for Basic Research
grant  99-01-0886}}
\author{Nikolai Dokuchaev\thanks{Correspondence address: 
The Institute of Mathematics and Mechanics,
St.Petersburg State University, Bibliotechnaya pl.2, Petrodvoretz,
St.Petersburg, 198904,  Russia. E-mail:
dokuchaev@pobox.spbu.ru; 
fax 7(812)4286998}\\
{\sm The Institute of Mathematics and Mechanics},
{\sm St.Petersburg State University, Russia}
}
\begin{document}
\maketitle
\begin{abstract}
We consider stochastic  
diffusion processes absorbed at the boundary of a domain. 
It is shown that there exist  initial distributions which ensure
a given  decreasing of density of the absorbed process.
\end{abstract}
\section{Introduction}
The paper investigates
 distributions of diffusion processes with adsorption 
on the boundary of a domain.  We consider
the problem of controlled absorption, and we obtain a new 
property of absorbed process: we establish an  existence of initial 
distributions  which ensure
a given  decreasing of density of the absorbed process.
More precisely,  we show that,
for any piecewise continuous function $\g(x)\ge 0$ such
that $\g(\cdot)\neq 0$,
there exists an initial probability density function  $p(\cdot,0)$
and a number $\a>0$ such that $p(x,0)\equiv p(x,T)+\a \g(x)$,
where $p(x,t)$ is
the probability density function of a process
$y(t)$ absorbed at the boundary $\p D$.
Thus, it can be concluded
that decreasing of density of the absorbed process can be
programmed.
\par
It can be mentioned that some problems for controlled
absorption were solved by author (1994, 1995) for
another type of constraints, when 
$v(x,T)\equiv \mu v(x,0)$ with $\mu\in (0,1)$,
where  $v(x,t)$ is a solution of a backward Kolmogorov equation.
 \section{Definitions}
Consider a $n$-dimensional Wiener process
$w(t)$ with independent components
on a complete probability space $(\Omega,{\cal F},\P)$,
such that $w(0)=0$.
\par
Let  $D\in \R^n$ be a bounded domain
with $C^2$ - smooth boundary $\p D$,
and let $T>0$ be a fixed number.
Consider the following It\^o's stochastic
differential equation:
\be
\label{y1}
\left\{\begin{array}{l}
dy(t)=f(y(t),t)dt+\b(y(t),t)dw(t), \quad t\in [0,T],\\
y(0)=y_0.
\end{array}\right.
\ee
Here $y_0$ is a random vector which is independent of
$w(\cdot)$ and such that $y\in D$ a.s.
\par
We assume that the function
$f(x,t):  \R^n \times  \R\to \R^{n}$ is measurable
and bounded, the function
$\b(x,t):  \R^n\times  \R\to \R^{n\times n}$
is continuous, and  that there exist
bounded derivatives
$\p^k \b(x,t)/\p x_i^k$,
where $i=1,...,n$,
 $k=1,2$.
We assume also that
$\b(x,t)\b(x,t)^\top\ge \d I_n$, where $\d>0$ is a constant,
and  $I_n$ is the unit matrix.
\par
Under these assumptions, there
exists the unique
weak solution  $y(t)$ 
of (\ref{y1}) (see e.g. Gihman and Skorohod (1979)).
We consider the process $y(t)$ which is absorbed at $\p D$,
i.e. until first exit from $D$.
\par\medskip
{\bf Spaces and classes of functions}\\
For a Banach space $X$,
we denote  the norm by $\|\cdot\|_{ X}$.
  \par
Let $H^0\defi L_2(D)$ and  $H^1\defi
\stackrel{0}{W_2^1}(D)$ be the standard Sobolev Hilbert spaces.
Let $H^{-1}$ be the dual space to $H^{1}$, with the norm
$\| \,\cdot\,\| _{H^{-1}}$ such that
$\| u\|_{ H^{-1}}$ for  \ $u \in H^{0}$
is supremum of $(u,v)_{H^0}$ over all
$v \in H^0$ such that $\| v\|_{H^1} \le 1 $;
we have $H^1\subset H^0\subset H^{-1}$, assuming the standard
embedding.
\par
We shall denote 
 the Lebesgue measure and
 the $\sigma $-algebra of Lebesgue sets in $\R ^m$ 
by $\oo\ell _{m}$ and $ {\oo{\cal B}}_{m}$, respectively.
\par
Introduce the spaces
$$
C^{k}(s,T)\defi C\left([s,T]; H^k\right),\quad
\W^{k}(s,T)\defi L^{2}\bigl([ s,T ],\oo{\cal B}_1,
\oo\ell_{1};  H^{k}\bigr),
\quad  k=0,\pm 1, 
$$
and the space 
$$
\V^{k}(s,T)\defi \W^{k}(s,T)\!\cap C^{k-1}(s,T), \quad k>0
$$
with the  norm 
$
\| u\| _{\V^k(s,T)} \defi \| u\| _{{\W}^k(s,T)} +\| u\| _{C^{k-1}(s,T)}.
$
\par
The following definition will be useful.
\begin{definition}
\label{def1}
A function $\g(\cdot):D\to\R$ is said to be piecewise continuous
if there exists a integer $N>0$ and a set of open 
domains $\{D_i\}_{i=1}^N$ such that
the following holds:
\begin{itemize}
\item
$D=\cup_{i=1}^ND_i$, $D_i\cap D_j=\emptyset $ for $i\neq j$;
\item for any $i\in\{1,...,N\}$, the function
$\g|_{D_i}$ is continuous and can be continued to a continuous function 
$\oo \g_i: D_i\cup\p D_i\to\R$;
\item for any $x\in \cup_{i=1}^N\p D_i$, there exists
$j\in\{1,...,N\}$ such that $x\in \p D_j$ and $\oo \g_j(x)=\g(x)$.
\end{itemize} 
\end{definition}
\section{The result} 
\begin{theorem}
\label{Thy}
Let $\g(\cdot):\R^n\to\R$ be a piecewise continuous function
such that $\g(x)\ge 0$
{\rm($\forall x$)} and $\g(\cdot)\neq 0$. Then
there exists
a random vector $y_0$ which is independent of
$w(\cdot)$ and  such that the following holds:
\par
{\rm (i)} the process $y(t)$,
defined by (\ref{y1}) and
absorbed at $\p D$,
has the probability density function
$p(x,t)\in \V^1(0,T)$; and
\par
{\rm (ii)} there exists a constant $\a>0$ such that
$p(x,0)\equiv p(x,T)+\a \g(x)$.
\end{theorem}
\par
{\it Proof}.
Let $y_0$ have a probability density function
$\rho(x)\in L_2(D)$. In that case, 
 the
probability density function $p(x,t)$ of the process $y(t)$
absorbed at $\p D$ satisfies the the forward Kolmogorov equation
\be
\label{y4}
\left\{
\begin{array}{ll}
\frac{\p p}{\p t}=A p, \quad t>0,\\
p(x,t)|_{x\in \p D}=0, \quad p(\cdot,0)=\rho(\cdot),
\end{array}
\right.
\ee
where
$$
Ap\defi\sum_{i,j=1}^n\frac{\p^2 }{\p x_i \p x_j}
\left(a_{ij}(x,t)p(x)\right)
-\sum_{i=1}^n\frac{\p}{\p x_i }\left( f_i(x,t)p(x)\right)
$$
and $a(x,t)\defi\{a_{ij}(x,t)\}=\frac{1}{2}\b(x,t)\b(x,t)^\top$.
\par
For  $s\in [0,T]$ and $\xi\in H^0$, consider the following auxiliary  
boundary value problem:
\be
\label{4.4}
\left\{
\begin{array}{ll}
\frac{dv}{dt}=A\, v, \quad t>s, \\
v(x,t)|_{x\in \p D}=0,
\quad  v(\cdot,s)=\xi(\cdot).
\end{array}
\right.
\ee
From the classic theory of parabolic equations,
it follows  that
\be
\label{1}
\|v(\cdot,T)\|_{H^1}\le C_1\|v(\cdot,s)\|_{H^1},
\ee
where $C_1>0$ is a constant which does not depend on $\xi$ and $s$
(see e.g. Ladyzenskaya {\it et al} (1968)).
\par
Introduce  operators $\L_s:H^0\to \V^1(s,T)$,
such that $\L_s\xi=v$,
where
$v$ is the solution in $\V^1(s,T)$ of the problem
(\ref{4.4}).
These linear operators are continuous 
(see e.g. Ladyzenskaya {\it et al} (1968)).
Introduce  an operator $\Q:H^{0}\to H^0$, such that
$\Q\xi=v(\cdot,T)$, where
$v=\L_0\xi$. Clearly, this operator is linear and continuous.
\begin{lemma}
\label{lemma1}
{\rm (i)} The operator
$\Q :H^0\to H^0$ is compact; \par
{\rm (ii)} If the equation $\Q \xi=\xi$ has
the only solution $\xi=0$ in $H^0$, then the operator 
$(I-\Q)^{-1}:H^0\to H^0$
is continuous.  
\end{lemma}
\par
{\it Proof of Lemma  \ref{lemma1}}.
Let $\xi\in H^0$ and  $v\defi\L_0 \xi$, i.e.
$v$ is the solution of the problem (\ref{4.4}).
We have that
$v|_{t\in[s,T]}=\L_sv(\cdot,s)$
for all $s\in [0,T]$, hence
$$
\begin{array}{ll}
\|v(\cdot,T)\|_{H^1}&\le C_1\inf_{t\in[0,T]}
\|v(\cdot,t)\|_{H^1}\\
&\le \frac{C_1}{\sqrt{T}} \left(\int_0^T\|v(\cdot,t)\|_{H^1}^2dt\right)^{1/2}
\le \frac{C_2}{\sqrt{T}}\|v\|_{\W^1}
\le \frac{C_3}{\sqrt{T}}\|\xi\|_{H^0}
\end{array}
$$
for constants $C_i>0$ which do not  depend on $\xi$.
Hence the operator $\Q:H^0\to H^1$ is continuous.
The embedding $H^1$ to $H^0$  is a compact operator
(see e.g. Yosida (1965), Ch.10.3). Then (i) follows.
Further, (ii) follows from Fredholm Theorem.
This completes the proof of Lemma \ref{lemma1}. $\Box$
\begin{lemma}
\label{lemma2}
For any $\g\in H^0$, there exists the unique
solution $u\in\V^1$ of the following problem:
\be
\label{y5}
\left\{
\begin{array}{ll}
\frac{\p u}{\p t}= A u, \quad t>0, \\
 u(x,t)|_{x\in \p D}=0,
\quad  u(x,0)-u(x,T)\equiv\g(x).
\end{array}
\right.
\ee
\end{lemma}
\par
{\it Proof of Lemma  \ref{lemma2}}.
First, we shall show that if $\g(\cdot)= 0$ then the unique solution
of (\ref{y5}) in $\V^1$ is $u(\cdot)=0$.
\par
Let   $u\in \V^1$ solve (\ref{y5}) with $\g(\cdot)= 0$.
Clearly, $u=\L_0u(\cdot,0)$. Denote $\zeta^+(x)\defi\max(0,u(x,0))$
and $\zeta^-(x)\defi\max(0,-u(x,0))$. 
Denote $u^-\defi\L_0\zeta^-$ and
$u^+\defi\L_0\zeta^+$.
We have that $u^+\ge 0$ and  $u^-\ge 0$ a.e.,
$u=u^+-u^-$ and  $u(x,0)\equiv\zeta^+(x)-\zeta^-(x)$.
\par
If  
$u(\cdot,0)\neq 0$ then either $\zeta^+(\cdot,0)\neq 0$ or 
$\zeta^-(\cdot,0)\neq 0$.
It follows from absorption at $\p D$
that if $\zeta^+(\cdot,0)\neq 0$ then 
$$
\int_Du^+(x,T)dx<\int_D\zeta^+(x)dx.
$$
(It suffices to note that the process $q\defi \nu^{-1} u^+$ is a
probability density function of a process with absorption on 
$\d D$, where $\nu\defi \int_D \zeta^+(x)dx$).  
Similarly,
if $\zeta^-(\cdot,0)\neq 0$ then 
$$
\int_Du^-(x,T)dx<\int_D\zeta^-(x)dx.
$$
Hence
$$
\begin{array}{ll}
\int_D|u(x,T)|dx &\le\int_Du^+(x,T)dx+\int_Du^-(x,T)dx\\
&<\int_D\zeta^+(x)dx+\int_D\zeta^-(x)dx=\int_D|u(x,0)|dx,
\end{array}
$$
i.e. 
$\int_D|u(x,T)|dx<\int_D|u(x,0)|dx$,
and the condition $u(x,0)\equiv u(x,T)$ fails to be  satisfied
for $u(\cdot)\neq 0$.
Thus, the unique solution
of (\ref{y5}) for $\g(\cdot)= 0$ is $u(\cdot)=0$.
\par
 By Lemma \ref{lemma1}, it follows that
the operator $(I-\Q)^{-1}:H^0\to H^0$ is continuous.
 Let $\g(\cdot)\in H^0$ be arbitrary.  
Then there exists $\zeta=(I-\Q)^{-1}\g\in H^0$,
and this $\zeta$ is unique. Let $u\defi\L_0\zeta$.
By the definitions of $\L_0$ and $\Q$, it follows 
that  $u(\cdot,T)=\Q u(\cdot,0)$.
We have that $u(\cdot,0)-\Q u(\cdot,0)=\g(\cdot)$, i.e. 
$u(\cdot,0)-u(\cdot,T)=\g(\cdot)$.
Thus, $u\defi\L_0\zeta=\L_0(I-\Q)^{-1}\g$ is 
the unique solution  
of (\ref{y5}) for any $\g(\cdot)\in H^0=L_2(D)$.
This completes the proof of Lemma \ref{lemma2}. $\Box$
\par
Let us continue the proof of the theorem.
Let $\g(x)\ge 0$ be such as in the assumptions of the theorem,
and let $u\defi\L_0(I-\Q)^{-1}\g$ solve the problem (\ref{y5}).
\par
 It is easy to see
that  if $u(\cdot,0)=0$ then $u(\cdot,T)=0$ and  $\g(\cdot)= 0$.
By the assumpltions,  
$\g(\cdot)\neq 0$, hence $u(\cdot,0)\neq 0$ and  $u(\cdot)\neq 0$.
\par
We remind that  $u=\L_0\zeta$, where $\zeta=u(\cdot,0)\in H^0$. 
By  Theorem 9.1
from Ch.IV of Ladyzenskaya {\it et al} (1968) applied for  smooth
functions which approximate $u(\cdot,0)$ in $H^0$,
and  by Theorem 8.1 from Ch.III
of this book,
 it follows that
there exists a  representative $\oo u(\cdot,T)$
of the corresponding element of $H^0$
which is continuous  in $x\in \oo D$, where $\oo D= D\cup\p D$.
(Note that, by the definition,  an element of $H^0=L_2(D)$  is a class of
$\oo\ell_n$-equivalent functions).
We have that  $u(\cdot,0)=u(\cdot,T)+\g(\cdot)$, hence there 
exists a  piecewise continuous representative $u'(\cdot,0)$
of $u(\cdot,0)\in H^0$. Note that 
such $u'(\cdot,0)$ 
 is not unique, because the boundary values
at $\p D_i$ 
in the Definition \ref{def1} can be choosen 
differently.
\par
Let us show that $u(x,0)\ge 0$ for a.e. $x$.
Suppose that
\be
\label{-}
\exists x\in D:\quad u'( x,0)<0.
\ee
If (\ref{-}) holds,  then there 
exists a 
piecewise continuous representative $\oo u(\cdot,0)$
such that
there exists $\w x\in D$ such that
$$
\oo u(\w x,0)<0, \quad
\oo u(\w x,0)\le \oo u( x,0) \quad \hbox{for a.e.}
\ x\in D.
$$
We have that
$$
\oo u(x,T)=\int_DG( x,y,T,0)u(y,0)dy,
$$
where $G(x,y,T,0)>0$ is the corresponding Green's function
for the problem (\ref{y1}).
Clearly,
$$
\G\defi\int_DG(\w x,y,T,0)dy\in (0,1).
$$
 Then
$$
\begin{array}{ll}
\oo u(\w x,T)&=\biggl(\int_D G(\w x,y,T,0)dy\biggr)\oo u(\w x,0)
+\int_DG(\w x,y,T,0)(\oo u(y,0)-\oo u(\w x,0))dy
\\
&\ge \G \oo u(\w x,0)>\oo u(\w x,0).
\end{array}
$$
It follows that if (\ref{-}) holds then
$u(x,t)$ does not satisfy  (\ref{y5}). Thus,
$u(x,0)\ge 0$ a.e.
\par
Now, let $\a\defi\left(\int_D u(x,0)dx\right)^{-1}$ and
$\rho(x)\defi\a \oo u(x,0)$.
We have that  $\rho(x)\ge 0$ and $\int_D\rho(x)dx=1$. Then
there exists a random vector   $y_0$ such that
$y_0$  
is independent of $w(\cdot)$ and has the
probability density function $\rho$.
Let $p\defi\L_0\rho$, i.e. $p$ solve (\ref{y4}). Clearly, 
$p(x,t)$ is the probability density function of the process 
$y(t)$ satisfying (\ref{y1}).  By linearity of (\ref{y5}), it follows 
that  $p=\a u$ and $p(\cdot,0)-p(\cdot,T)=\a \g(\cdot)$. 
This completes the proof of Theorem \ref{Thy}. $\Box$
\\
\
\\
{\Large\bf References}
 \par
Dokuchaev, N., 1994.
Parabolic  equations  without  Cauchy conditions  and
control problems for diffusion processes. Part I. {\it
Differential Equations} {\bf 30},  1606-1617.  
\par
  Dokuchaev, N., 1995.
Parabolic  equations  without  Cauchy conditions  and
control problems for diffusion processes. Part II. {\it
Differential Equations} 
{\bf  31},   1362-1372.  
\par
Gihman I.I. and Skorohod, A.V. 1979.
{\it The Theory of Stochastic Processes.} III. 
Springer-Verlag.
\par
Ladyzenskaya, O.A.,  Solonnikov, V.A., and   Ural'ceva, N.N., 1968.
{\it Linear and Quasi--Linear Equations of Parabolic Type.}
Providence, R.I.: American Mathematical Society.
\par
Yosida, K., 1965. {\it Functional Analysis}. Springer-Verlag. Berlin, Gottingen,
Heidelberg.
\\
\
\\

\end{document}